\documentclass[english]{article}
\usepackage[T1]{fontenc}
\usepackage[latin9]{inputenc}
\usepackage[a4paper]{geometry}
\usepackage{url}
\usepackage{amsmath}
\usepackage{amssymb}
\usepackage{graphicx}
\usepackage{setspace}
\usepackage{cite}
\usepackage{hyperref}
\usepackage{babel}
\DeclareMathOperator{\rmd}{d}
\DeclareMathOperator{\real}{Re}
\newcommand{\iu}{\mathrm{i}}

\begin{document}
\title{$1/f$ noise in semiconductors arising from the heterogeneous detrapping
process of individual charge carriers}
\author{Aleksejus Kononovicius\thanks{email: \protect\href{mailto:aleksejus.kononovicius@tfai.vu.lt}{aleksejus.kononovicius@tfai.vu.lt};
website: \protect\url{https://kononovicius.lt}} and Bronislovas Kaulakys}
\date{Institute of Theoretical Physics and Astronomy, Vilnius University}
\maketitle
\begin{abstract}
We propose a model of $1/f$ noise in semiconductors based on the
drift of individual charge carriers and their interaction with the
trapping centers. We assume that the trapping centers are homogeneously
distributed in the material. The trapping centers are assumed to be
heterogeneous and have unique detrapping rates. We show that uniform
detrapping rate distribution emerges as a natural consequence of the
vacant trap depths following the Boltzmann distribution, and the detrapping
process obeying Arrhenius law. When these laws apply, and if the trapping
rate is low in comparison to the maximum detrapping rate, $1/f$ noise
in the form of Hooge's relation is recovered. Hooge's parameter, $\alpha_{H}$,
is shown to be a ratio between the characteristic trapping rate and
the maximum detrapping rate. The proposed model implies that $1/f$
noise arises from the temporal charge carrier number fluctuations,
not from the spatial mobility fluctuations.
\end{abstract}

\section{Introduction}

The nature of the $1/f$ noise (often also referred to as low frequency,
flicker or pink noise), characterized by power spectral density of
$S\left(f\right)\sim1/f^{\beta}$ form (with $0.5\leq\beta\leq1.5$),
remains open to discussion despite almost 100 years since the first
reports \cite{Johnson1925PR,Schottky1926PR,Milotti2002,Kogan1996CUP}.
While many materials, devices, and systems exhibit different kinds
of fluctuations or noise \cite{Kogan1996CUP,Lowen2005Wiley,VanKampen2007NorthHolland},
only the white noise and the Brownian noise are well understood from
the first principles. White noise is characterized by absence of any
temporal correlations, and has a flat power spectral density of $S\left(f\right)\sim1/f^{0}$
form. Examples of the white noise include thermal and shot noise.
Thermal noise is known to arise from the random motion of the charge
carriers. It occurs at any finite temperature regardless of whether
the current flows. Shot noise, on the other hand, is a result of the
discrete nature of the charge carriers and the Poisson statistics
of waiting times before each individual detection of the charge carrier.
The Brownian noise is a temporal integral of the white noise, and
thus exhibits no correlations between the increments of the signal,
it is characterized by a power spectral density of $S\left(f\right)\sim1/f^{2}$
form.

Theory of $1/f$ noise based on the first principles is still an open
problem. $1/f$ noise is of particular interest as it is observed
across various physical \cite{Voss1976PRL,Dutta1981RMP,Mitin2002,Balandin2013NN,Fox2021Nature,Wirth2021IEEE},
and non-physical \cite{Voss1975Nature,Kobayashi1982BioMed,Gilden1995Science,Levitin2012PNAS}
systems. $1/f$ noise cannot be obtained by the simple procedure of
integration, differentiation, or simple transformations of well-understood
processes. Also the general mechanism of generating $1/f$ noise has
not yet been properly identified, and there is no generally accepted
solution to the $1/f$ noise problem.

The oldest explanation for $1/f$ noise involves the superposition
of Lorentzian spectra \cite{Bernamont1937ProcPhysSoc,Surdin1951,McWhorter1957,VanDerZiel1979AEEP}.
Lorentzian spectral densities themselves may arise from the random
telegraph signals \cite{Kogan1996CUP}, and from the Brownian motion
with a broad distribution of relaxations \cite{Kaulakys2005PhysRevE}.
These approaches, as well as many others, are often limited to the
specific systems being modeled, or require quite restrictive assumptions
to be satisfied \cite{Wong2003MR}. In the recent decades, series
of models for the $1/f$ noise based on the specific, autoregressive
$\mathrm{AR}\left(1\right)$, point process \cite{Kaulakys2005PhysRevE},
and the agent-based model \cite{Ruseckas2011EPL,Kononovicius2012PhysA},
yielding nonlinear stochastic differential equation \cite{Kaulakys2009JStatMech}
was proposed (see \cite{Kazakevicius2021Entropy} for a recent review).
Another more recent trend relies on scaling properties and nonlinear
transformations of signals \cite{Ruseckas2014JStatMech,Kaulakys2015MPLB,Eliazar2021JPA,Kazakevicius2023PRE}.
These models, on the other hand, prove to be rather more abstract,
and therefore more similar to the long-range memory models found in
the mathematical literature, such as fractional Brownian motion \cite{Mandelbrot2013Springer,Beran2017Routledge}
or ARCH models \cite{Giraitis2009,Bollerslev2023JEconom}. These and
other similar models of $1/f$ noise are hardly applicable to the
description and explanation of the mostly observable $1/f$ noise
in the semiconductors.

On the other hand, for a homogeneous semiconductor material Hooge
proposed an empirical relation for the $1/f$ noise dependence on
the parameters of the material \cite{Hooge1969PLA,Hooge1972Phys},
\begin{equation}
S\left(f\right)=\bar{I}^{2}\frac{\alpha_{H}}{Nf}.\label{eq:hooge-relation}
\end{equation}
Where $\bar{I}$ stands for the average current flowing through the
cross-section of the semiconductor material, $N$ is the number of
charge carriers, and $\alpha_{H}$ is the titular Hooge parameter.
If the current is kept constant, or does not exhibit large fluctuations,
Hooge's empirical relation could be rewritten in terms of voltage
or resistivity noise, i.e., $S_{V}\left(f\right)=\bar{V}^{2}\frac{\alpha_{H}}{Nf}$
or $S_{R}\left(f\right)=\bar{R}^{2}\frac{\alpha_{H}}{Nf}$ (here the
subscripts emphasize fluctuations of which quantity are being observed).
However, we are specifically interested in the case of constant voltage,
focusing on the power spectral density of the current fluctuations
that are associated with Eq.~(\ref{eq:hooge-relation}). There were
numerous attempts to derive or explain the structure of the Hooge's
relation \cite{Hooge1994IEEE,Kaulakys1999PLA,Dmitriev2009JAppPhys,Vandamme2013ICNF,Palenskis2015SRep}.
A more recent derivation of the Hooge's parameter, based on the Poisson
generation-recombination process modulated by random telegraph noise,
was conducted in \cite{Gruneis2019PLA,Gruneis2022PhysA}. Yet these
models, as well as many others, cannot be directly applied to describe
and explain the widespread $1/f$ noise in the semiconductors.

Here, we propose a model of $1/f$ noise in semiconductors containing
heterogeneous trapping centers. As far as the square of the average
current $\bar{I}^{2}$ is proportional to the squared number of the
charge carriers $N^{2}$, Hooge's relation implies that the intensity
of $1/f$ noise is proportional to the number of charge carriers $N$.
Therefore, as the first approximation we can consider the noise originating
from the flow of individual charge carriers. It is known that the
drift, and the diffusion, of the charge carriers does not yield $1/f$
noise \cite{Kogan1996CUP}. Therefore, we consider the drift of the
charge carriers interrupted by their entrapment in the trapping centers.
We show that, if the detrapping rates of individual trapping centers
are heterogeneous and uniformly distributed, $1/f$ noise arises.
As an explanation for the uniform detrapping rate distribution, we
note that it may arise from the interplay between the Boltzmann distribution
of the vacant trap depths (as is observed in various materials \cite{Bisquert2008PRB,Wong2020ACSEnergyLett,Beckers2021JAppPhys,Zeiske2021NatureComm})
and the Arrhenius law (which is often applied in empirical works studying
varied activation and detrapping processes in semiconductors \cite{Peters2015JPCB,Kurpiers2018NatComm,Arakawa2020NatMater,Kumar2024PRE}).
In this model, the signal generated by a single charge carrier is
similar to the signal composed of non-overlapping rectangular pulses
\cite{Kononovicius2023PRE}. Here, we derive Hooge's relation, and
show that Hooge's parameter is a ratio between the characteristic
trapping rate and the maximum detrapping rate. The proportionality
between Hooge\textquoteright s parameter and the characteristic trapping
rate was reported earlier in quite a few experimental works \cite{Lukyanchikova1990PhysB,Tacano2004FluctuationsNoiseMater,Tousek2024SciRep}.
This result prompts us to suggest that $1/f$ noise in semiconductors
arises from the fluctuations in the effective number of charge carriers,
not from the spatial fluctuations in mobility.

This paper is organized as follows. In Section~\ref{sec:general-model}
we introduce a model for $1/f$ noise in the semiconductors based
on the trapping-detrapping process of a single charge carrier. In
Section~\ref{sec:finite-experiments} we address the implications
of finite experiments and simulations. Namely, we show that the power
spectral density produced by a single charge carrier may exhibit spurious
low-frequency cutoff. This cutoff disappears, if the current generated
by a large number of charge carriers is considered. Finally, Hooge's
empirical relation and Hooge's parameter value for the proposed model
is derived in Section~\ref{sec:hooge-expression}. The main results
of the paper are summarized in Section~\ref{sec:conclusions}.

\section{Model for $1/f$ noise in a homogeneous semiconductor material\label{sec:general-model}}

\begin{figure}[h]
\begin{centering}
\includegraphics[width=0.9\textwidth]{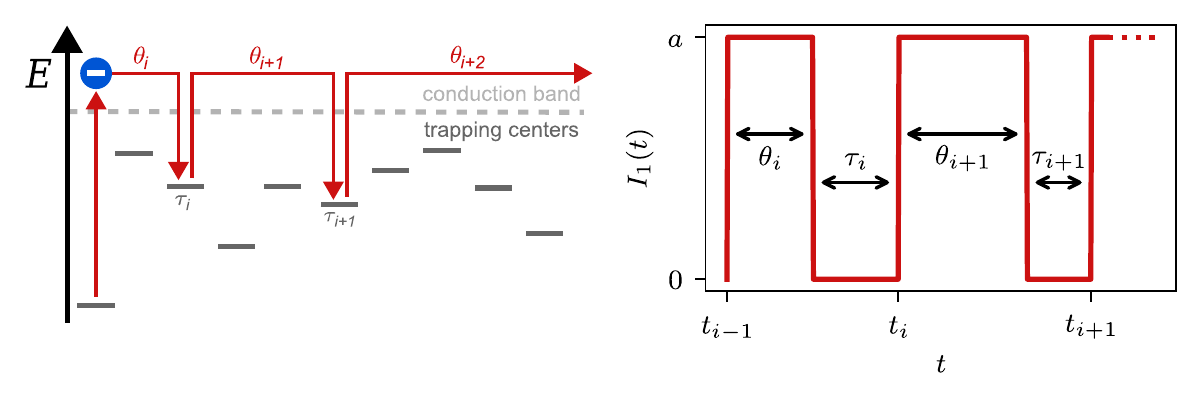}
\par\end{centering}
\caption{Visualization of the single charge carrier trapping-detrapping process
(left) and a sample single charge carrier contribution to the net
current (right). Relevant notation: $\tau_{i}$ is the detrapping
time (gap duration), $\theta_{i}$ is the trapping time (pulse duration),
$a$ is the height of the pulses (single free charge carrier contribution
to the net current), $t_{i}$ is the time of $i$-th detrapping event.\label{fig:explanation}}
\end{figure}

Let us consider a drift of a single charge carrier (e.g., election)
through a homogeneous semiconductor material. While the charge carrier
is freely moving through the conduction band, it will generate a non-zero
contribution to the net current, i.e., $I_{1}\left(t\right)=a$ for
$t$ when the charge carrier is free. As the material contains trapping
centers, the freely moving charge carrier will eventually get trapped
in one of such trapping centers. Let $\tau_{i}$ stand for $i$-th
detrapping time (time spent in the trap) and $\theta_{i}$ be $i$-th
trapping time (time spent moving). Under these considerations the
contribution of single charge carrier to the net current will be composed
of gaps (duration corresponds to the respective detrapping time) and
pulses (duration corresponds to the respective trapping time). For
visual illustration of the single charge carrier trapping-detrapping
process and a sample signal see Fig.~\ref{fig:explanation}.

The power spectral density of a signal with rectangular pulses of
fixed height is given by \cite{Kononovicius2023PRE}
\begin{equation}
S_{1}\left(f\right)=\lim_{T\rightarrow\infty}\left\langle \frac{2}{T}\left|\int_{0}^{T}I_{1}\left(t\right)e^{-2\pi\iu ft}\rmd t\right|\right\rangle =\frac{a^{2}\bar{\nu}}{\pi^{2}f^{2}}\real\left[\frac{\left(1-\chi_{\theta}\left(f\right)\right)\left(1-\chi_{\tau}\left(f\right)\right)}{1-\chi_{\theta}\left(f\right)\chi_{\tau}\left(f\right)}\right].\label{eq:rectangular-psd}
\end{equation}
In the above $T$ stands for observation time (duration of the signal),
which is assumed to approach infinity \cite{Kononovicius2023PRE},
$\chi_{\tau}\left(f\right)$ and $\chi_{\theta}\left(f\right)$ stand
for the characteristic functions of the respective detrapping and
trapping time distributions, while $\bar{\nu}$ is the mean number
of pulses per unit time. For the ergodic processes, and given a long
observation time $T$, the value of $\bar{\nu}$ is trivially derived
from the mean trapping and detrapping times, i.e., $\bar{\nu}=\frac{1}{\left\langle \theta\right\rangle +\left\langle \tau\right\rangle }$.
For the nonergodic processes, or if the observation time $T$ is comparatively
short, the expected value of $\bar{\nu}$ can be derived from the
means of the appropriately truncated distributions, or it may be defined
purely empirically, i.e., $\bar{\nu}=K/T$ (here $K$ is the number
of observed pulses).

Typically when trapping-detrapping processes are considered \cite{Kogan1996CUP,Mitin2002}
it is assumed that both $\tau_{i}$ and $\theta_{i}$ are sampled
from the exponential distributions with rates $\gamma_{\tau}$ and
$\gamma_{\theta}$ respectively. Characteristic function of the exponential
distribution with an event rate $\gamma$, is given by
\begin{equation}
\chi\left(f\right)=\int_{0}^{\infty}\gamma e^{2\pi\iu f\tau-\gamma\tau}\rmd\tau=\frac{\gamma}{\gamma-2\pi\iu f}.\label{eq:exp-characteristic}
\end{equation}
Inserting Eq.~(\ref{eq:exp-characteristic}) as the characteristic
function for both trapping and detrapping time distributions into
Eq.~(\ref{eq:rectangular-psd}) yields a Lorentzian power spectral
density \cite{Kogan1996CUP}. Notably, there were prior works which
have examined the case when $\tau_{i}$, $\theta_{i}$, or both are
sampled from distributions with power-law tails \cite{Margolin2006JStatPhys,Lukovic2008JChemPhys,Niemann2013PRL,Leibovich2016PRE,Gruneis2019PLA,Gruneis2022PhysA,Kononovicius2023PRE}.
Under the power-law distribution assumption, it was shown $S\left(f\right)\sim1/f^{\beta}$
dependence can be recovered.

Here, let us assume that the trapping centers are heterogeneous. Each
of them has their own unique depth, or detrapping (activation) energy,
$E_{a}^{\left(i\right)}$. As is commonly observed \cite{Peters2015JPCB,Kurpiers2018NatComm,Arakawa2020NatMater,Kumar2024PRE},
let us assume that the detrapping process obeys Arrhenius law
\begin{equation}
\gamma_{\tau}^{\left(i\right)}=A\exp\left[-\frac{E_{a}^{\left(i\right)}}{k_{B}\Theta}\right].
\end{equation}
To obtain the overall detrapping time distribution we first need to
establish the distribution of detrapping energies. Not all trapping
centers will participate in the trapping-detrapping process at all
times. Because charge carrier first needs to be trapped, before being
detrapped, only vacant trapping centers will participate in the process.
In experimental literature \cite{Bisquert2008PRB,Wong2020ACSEnergyLett,Beckers2021JAppPhys,Zeiske2021NatureComm}
it is well established that vacant trap level depths (their activation
energies) reasonably well follow the Boltzmann distribution 
\begin{equation}
p\left[E_{a}^{\left(i\right)}\right]=C_{N}\exp\left[-\frac{E_{a}^{\left(i\right)}}{k_{B}\Theta}\right].
\end{equation}
In the above $C_{N}$ stands for the normalization constant. Notably,
this result also follows directly from the Fermi-Dirac statistics
under the assumption that trap level degeneracy is constant in respect
to activation energy. Then, from the conservation of the probability
density, it follows that the distribution of detrapping rates would
be uniform
\begin{equation}
p\left[\gamma_{\tau}^{\left(i\right)}\right]=\frac{p\left[E_{a}^{\left(i\right)}\right]}{\left|\frac{\rmd c}{\rmd E_{a}^{\left(i\right)}}\right|}=\frac{C_{N}\exp\left[-\frac{E_{a}^{\left(i\right)}}{k_{B}\Theta}\right]}{\frac{A}{k_{B}\Theta}\exp\left[-\frac{E_{a}^{\left(i\right)}}{k_{B}\Theta}\right]}=\mathrm{const}.
\end{equation}
It is important to note that other physical mechanisms could also
imply uniform distribution of the detrapping rates as long as $\frac{p\left(\eta\right)}{\left|\frac{\rmd\gamma_{\tau}^{\left(i\right)}}{\rmd\eta}\right|}=\mathrm{const}$
(here $\eta$ is some generic physical quantity which would impact
the detrapping process).

\begin{figure}[ht]
\begin{centering}
\includegraphics[width=0.45\textwidth]{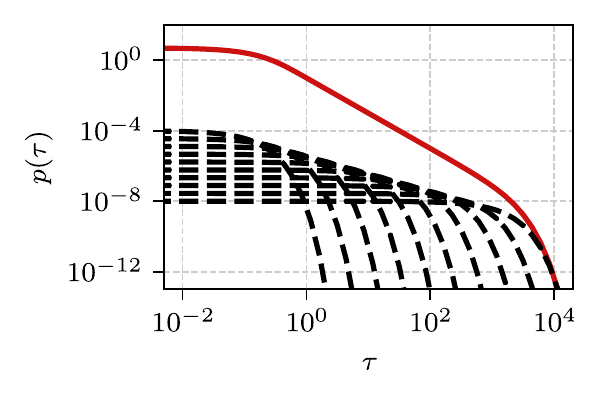}
\par\end{centering}
\caption{Probability density function of the detrapping time distribution under
the assumption that detrapping rates of individual trapping centers
are uniformly distributed (red curve), Eq.~(\ref{eq:escape-time-pdf}).
The probability density function was calculated for $\gamma_{\text{min}}=10^{-3}$,
and $\gamma_{\text{max}}=10$ case. Black dashed curves correspond
to the exponential probability density functions of the detrapping
times from the individual trapping centers with fixed rates: $\gamma_{\tau}=10^{-3}$,
$2.78\times10^{-3}$, $7.74\times10^{-3}$, $2.15\times10^{-2}$,
$5.99\times10^{-2}$, $1.67\times10^{-1}$, $4.64\times10^{-1}$,
$1.29$, $3.59$, and $10$. Normalization of the exponential probability
density functions was adjusted for the visualization purposes, but
it remains proportional to their respective contributions.\label{fig:distribution}}
\end{figure}

Let $\gamma_{\tau}^{\left(i\right)}$ be uniformly distributed in
$\left[\gamma_{\text{min}},\gamma_{\text{max}}\right]$. Then it can
be shown that the probability density function of the detrapping time
distribution is given by
\begin{equation}
p\left(\tau\right)=\frac{1}{\gamma_{\text{max}}-\gamma_{\text{min}}}\int_{\gamma_{\text{min}}}^{\gamma_{\text{max}}}\gamma_{\tau}\exp\left(-\gamma_{\tau}\tau\right)\rmd\gamma_{\tau}=\frac{\left(1+\gamma_{\text{min}}\tau\right)\exp\left(-\gamma_{\text{min}}\tau\right)-\left(1+\gamma_{\text{max}}\tau\right)\exp\left(-\gamma_{\text{max}}\tau\right)}{\left(\gamma_{\text{max}}-\gamma_{\text{min}}\right)\tau^{2}}.\label{eq:escape-time-pdf}
\end{equation}
This probability density function saturates for the short detrapping
times, $\tau\ll\frac{1}{\gamma_{\text{max}}}$. For the longer detrapping
times, $\tau\gg\frac{1}{\gamma_{\text{min}}}$, it decays as an exponential
function. In the intermediate value range, $\frac{1}{\gamma_{\text{max}}}\ll\tau\ll\frac{1}{\gamma_{\text{min}}}$,
this probability density function has the $\tau^{-2}$ asymptotic
behavior, which is already known to lead to $1/f$ noise \cite{Margolin2006JStatPhys,Lukovic2008JChemPhys,Niemann2013PRL,Kononovicius2023PRE}.
The benefit of this formulation is that it allows to see how the $\tau^{-2}$
asymptotic behavior can emerge in homogeneous semiconductors. Experimentally
$\tau^{-2}$ asymptotic behavior is observable in quantum dots, nanocrystal,
nanorod, and other semiconductors \cite{Frantsuzov2008NatPhys,Cordones2013CSR,Nenashev2018PRB,Haneef2020JMCC},
with the detrapping times ranging from picoseconds to several months.
The asymptotic behavior of Eq.~(\ref{eq:escape-time-pdf}) can be
examined in Fig.~\ref{fig:distribution} where it is represented
by a red curve. Fig.~\ref{fig:distribution} also highlights contributions
of some of the individual trapping centers, detrapping time distributions
of which are plotted as dashed black curves.

Unlike the simple power-law distribution, this detrapping time distribution
does not require the introduction of any arbitrary cutoffs. Also the
parameters of this detrapping time distribution have explicit physical
meaning. Furthermore, the statistical moments are well-defined and
have compact analytical forms. The mean of the distribution is given
by
\begin{equation}
\left\langle \tau\right\rangle =\frac{1}{\gamma_{\text{max}}-\gamma_{\text{min}}}\ln\left(\frac{\gamma_{\text{max}}}{\gamma_{\text{min}}}\right).\label{eq:mean-gap}
\end{equation}
Higher order moments also exist and can be easily derived.

The characteristic function of the detrapping time distribution can
be obtained either by calculating Fourier transform of Eq.~(\ref{eq:escape-time-pdf}),
or by averaging over the characteristic functions of the exponential
distribution, Eq.~(\ref{eq:exp-characteristic}). Both approaches
lead to the same expression, but the latter approach is quicker
\begin{equation}
\chi_{\tau}\left(f\right)=\frac{1}{\gamma_{\text{max}}-\gamma_{\text{min}}}\int_{\gamma_{\text{min}}}^{\gamma_{\text{max}}}\frac{\gamma_{\tau}}{\gamma_{\tau}-2\pi\iu f}\rmd\gamma_{\tau}=1+\frac{2\pi\iu f}{\gamma_{\text{max}}-\gamma_{\text{min}}}\ln\left(\frac{\gamma_{\text{max}}-2\pi\iu f}{\gamma_{\text{min}}-2\pi\iu f}\right).
\end{equation}
If the interval of the possible detrapping rates is broad $\gamma_{\text{min}}\ll\gamma_{\text{max}}$,
then for $\gamma_{\text{min}}\ll2\pi f\ll\gamma_{\text{max}}$ the
characteristic function can be approximated by
\begin{equation}
\chi_{\tau}\left(f\right)\approx1+\frac{2\pi\iu f}{\gamma_{\text{max}}}\ln\left(1+\frac{\iu\gamma_{\text{max}}}{2\pi f}\right)\approx1-\frac{2\pi f}{\gamma_{\text{max}}}\left[\frac{\pi}{2}-\iu\ln\left(\frac{2\pi f}{\gamma_{\text{max}}}\right)\right].\label{eq:uniform-escape-char}
\end{equation}
Inserting Eq.~(\ref{eq:uniform-escape-char}) into Eq.~(\ref{eq:rectangular-psd})
we have
\begin{equation}
S_{1}\left(f\right)=\frac{2a^{2}\bar{\nu}}{\pi\gamma_{\text{max}}f}\real\left[\frac{\left(1-\chi_{\theta}\left(f\right)\right)\left[\frac{\pi}{2}-\iu\ln\left(\frac{2\pi f}{\gamma_{\text{max}}}\right)\right]}{1-\chi_{\theta}\left(f\right)\left\{ 1-\frac{2\pi f}{\gamma_{\text{max}}}\left[\frac{\pi}{2}-\iu\ln\left(\frac{2\pi f}{\gamma_{\text{max}}}\right)\right]\right\} }\right].
\end{equation}
Assuming that $\frac{2\pi f}{\gamma_{\text{max}}}\left[\frac{\pi}{2}-\iu\ln\left(\frac{2\pi f}{\gamma_{\text{max}}}\right)\right]\ll1$,
which is supported by an earlier assumption that $2\pi f\ll\gamma_{\text{max}}$,
allows to simplify the above to
\begin{equation}
S_{1}\left(f\right)\approx\frac{a^{2}\bar{\nu}}{\gamma_{\text{max}}f}.\label{eq:main-result}
\end{equation}
This approximation should hold well for $\gamma_{\text{min}}\ll2\pi f\ll\gamma_{\text{max}}$,
and should not depend on the explicit form of $\chi_{\theta}\left(f\right)$
unless $\chi_{\theta}\left(f\right)\approx1$ for at least some of
the frequencies in the range.

\begin{figure}[ht]
\begin{centering}
\includegraphics[width=0.45\textwidth]{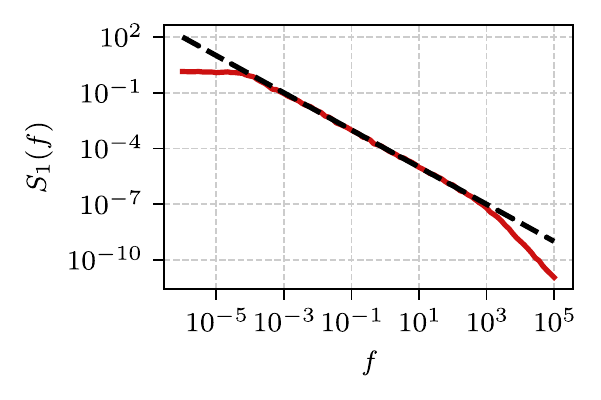}
\par\end{centering}
\caption{Power spectral density of the simulated signal (red curve) and its
analytical approximation by Eq.~(\ref{eq:main-result}) (black dashed
curve). Simulated power spectral density was obtained by averaging
over $10^{2}$ realizations. Simulation parameters: $T=10^{6}$, $\gamma_{\text{min}}=10^{-4}$,
$\gamma_{\text{max}}=10^{4}$, $a=1$, $\gamma_{\theta}=1$.\label{fig:sample-psd}}
\end{figure}

Let us examine a specific case when the trapping centers are uniformly
distributed within the material, and therefore the trapping process
can be assumed to be a homogeneous Poisson process. Inserting the
characteristic function of the exponential distribution, Eq.~(\ref{eq:exp-characteristic}),
as the characteristic function of the trapping time distribution into
Eq.~(\ref{eq:rectangular-psd}) yields
\begin{equation}
S_{1}\left(f\right)=\frac{4a^{2}\bar{\nu}}{\gamma_{\theta}^{2}}\real\left[\frac{1}{1-\chi_{\tau}\left(f\right)-\frac{2\pi\iu f}{\gamma_{\theta}}}\right].\label{eq:poisson-pulse-psd}
\end{equation}
Then inserting the characteristic function of the proposed detrapping
time distribution, Eq.~(\ref{eq:uniform-escape-char}), into Eq.~(\ref{eq:poisson-pulse-psd})
yields
\begin{equation}
S_{1}\left(f\right)=\frac{a^{2}\bar{\nu}\gamma_{\text{max}}}{\gamma_{\theta}^{2}f}\times\frac{1}{\left(\frac{\pi}{2}\right)^{2}+\left[\frac{\gamma_{\text{max}}}{\gamma_{\theta}}+\ln\left(\frac{2\pi f}{\gamma_{\text{max}}}\right)\right]^{2}}.\label{eq:exp-main-result}
\end{equation}
If the maximum detrapping rate is large in comparison to the trapping
rate, i.e., $\frac{\gamma_{\text{max}}}{\gamma_{\theta}}\gg\frac{\pi}{2}$
and $\frac{\gamma_{\text{max}}}{\gamma_{\theta}}\gg-\ln\left(\frac{2\pi f}{\gamma_{\text{max}}}\right)$,
then we recover Eq.~(\ref{eq:main-result}). In Fig.~\ref{fig:sample-psd}
the power spectral density of a simulated signal with comparatively
large detrapping rates is shown as a red curve. We have chosen observation
time $T$ to allow us to show three regimes of the power spectral
density: white noise cutoff for $2\pi f\ll\gamma_{\text{min}}$, $1/f$
noise for $\gamma_{\text{min}}\ll2\pi f\ll\gamma_{\text{max}}$ and
Brown noise for $\gamma_{\text{max}}\ll2\pi f$. Longer or similar
observation times would yield similar power spectral density.

\section{Low-frequency cutoff in finite experiments\label{sec:finite-experiments}}

The obtained approximation, Eq.~(\ref{eq:main-result}), holds in
the infinite observation time limit (single signal of infinite duration
$T$) or the infinite number of experiments limit (infinitely many
signals with finitely long observation time $T$). If either of the
limits doesn't hold, then the range of frequencies over which the
pure $1/f$ noise is observed becomes narrower. In the finite experiments
the process will not reach a steady state, and therefore the cutoff
frequencies will depend not on the model parameter values $\gamma_{\text{min}}$
and $\gamma_{\text{max}}$, but on the smallest and the largest $\gamma_{\tau}^{\left(i\right)}$
values actually observed during the experiment. The difference between
$\gamma_{\text{max}}$ and the largest $\gamma_{\tau}^{\left(i\right)}$
is negligible, because the pure $1/f$ noise will be observed only
if $\gamma_{\text{max}}$ is a relatively large number. On the other
hand the relative difference between $\gamma_{\text{min}}$ and smallest
$\gamma_{\tau}^{\left(i\right)}$ might not be negligible. Let us
estimate the expected value of the smallest $\gamma_{\tau}^{\left(i\right)}$
in a finite experiment.

\begin{figure}[ht]
\begin{centering}
\includegraphics[width=0.45\textwidth]{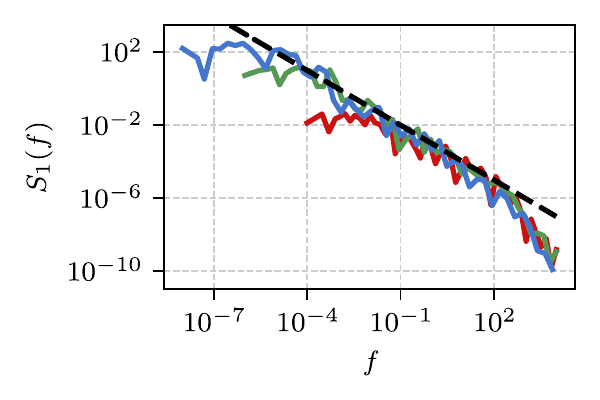}
\par\end{centering}
\caption{The effect of increasing the observation time $T$
on the obtained power spectral density. Dashed black curve corresponds
to Eq.~(\ref{eq:main-result}). Simulation parameters: $a=1$, $\gamma_{\theta}=1$,
$\gamma_{\text{min}}=0$, $\gamma_{\text{max}}=10^{3}$, $T=10^{4}$
(red curve), $10^{6}$ (green curve), and $10^{8}$ (blue curve).\label{fig:different-duration}}
\end{figure}

In the model introduced in the previous section $\gamma_{\tau}^{\left(i\right)}$
is sampled from the uniform distribution with $\left[\gamma_{\text{min}},\gamma_{\text{max}}\right]$
range of possible values. It is known that, for $x_{i}$ sampled from
the uniform distribution with $\left[0,1\right]$ range of possible
values, the smallest $x_{i}$ observed in the sample of size $K$
is distributed according to the Beta distribution with the shape parameters
$\alpha_{1}=1$ and $\alpha_{2}=K$ \cite{Gentle2009Springer}. Thus
the expected value of the smallest $x_{i}$ is given by
\begin{equation}
\left\langle \min\left\{ x_{i}\right\} _{K}\right\rangle =\frac{\alpha_{1}}{\alpha_{1}+\alpha_{2}}=\frac{1}{K+1}.
\end{equation}
Rescaling the range of possible values to $\left[\gamma_{\text{min}},\gamma_{\text{max}}\right]$
yields
\begin{equation}
\gamma_{\text{min}}^{\text{(eff)}}=\left\langle \min\left\{ \gamma_{\tau}^{\left(i\right)}\right\} _{K}\right\rangle =\frac{\gamma_{\text{max}}-\gamma_{\text{min}}}{K+1}+\gamma_{\text{min}}.
\end{equation}
As $K$ corresponds to the number of pulses in the signal, we have
that $K=\bar{\nu}T=\frac{T}{\left\langle \theta\right\rangle +\left\langle \tau\right\rangle }$
and
\begin{equation}
\gamma_{\text{min}}^{\text{(eff)}}=\left(\gamma_{\text{max}}-\gamma_{\text{min}}\right)\frac{\left\langle \theta\right\rangle +\left\langle \tau\right\rangle }{\left\langle \theta\right\rangle +\left\langle \tau\right\rangle +T}+\gamma_{\text{min}}.
\end{equation}
In the above $\left\langle \theta\right\rangle $ is effectively a
model parameter as it is trivially given by $\left\langle \theta\right\rangle =\frac{1}{\gamma_{\theta}}$,
while $\left\langle \tau\right\rangle $ is a derived quantity which
has a more complicated dependence on the model parameters $\gamma_{\text{min}}$
and $\gamma_{\text{max}}$ (see Eq.~(\ref{eq:mean-gap})). If the
range of possible $\gamma_{\tau}^{\left(i\right)}$ values is broad,
i.e., $\gamma_{\text{max}}\gg\gamma_{\text{min}}$, we have
\begin{equation}
\gamma_{\text{min}}^{\text{(eff)}}\approx\gamma_{\text{max}}\frac{\gamma_{\text{max}}\left\langle \theta\right\rangle +\ln\frac{\gamma_{\text{max}}}{\gamma_{\text{min}}}}{\gamma_{\text{max}}\left(\left\langle \theta\right\rangle +T\right)+\ln\frac{\gamma_{\text{max}}}{\gamma_{\text{min}}}}+\gamma_{\text{min}}.
\end{equation}
The above applies to the ergodic case with $\gamma_{\text{min}}\gg1/T$.
In the nonergodic case, for $\gamma_{\text{min}}\lesssim1/T$, it
would impossible to distinguish between the cases corresponding to
the different $\gamma_{\text{min}}$ values. Therefore, for the nonergodic
case, $\gamma_{\text{min}}$ can be replaced by $1/T$ yielding
\begin{equation}
\gamma_{\text{min}}^{\text{(eff)}}\approx\gamma_{\text{max}}\frac{\gamma_{\text{max}}\left\langle \theta\right\rangle +\ln\left(\gamma_{\text{max}}T\right)}{\gamma_{\text{max}}\left(\left\langle \theta\right\rangle +T\right)+\ln\left(\gamma_{\text{max}}T\right)}+\frac{1}{T}\approx\frac{1+\gamma_{\text{max}}\left\langle \theta\right\rangle +\ln\left(\gamma_{\text{max}}T\right)}{T}.
\end{equation}
For relatively long trapping times, $\left\langle \theta\right\rangle \gg\frac{\ln\left(\gamma_{\text{max}}T\right)}{\gamma_{\text{max}}}$,
we have that
\begin{equation}
\gamma_{\text{min}}^{\text{(eff)}}\approx\frac{1+\gamma_{\text{max}}\left\langle \theta\right\rangle }{T}\approx\frac{\gamma_{\text{max}}}{\gamma_{\theta}T}.
\end{equation}
From the above, it follows that low-frequency cutoff is always present
in singular experiments with one charge carrier, and with finite observation
time $T$. The cutoff will be observed at a frequency close to $\gamma_{\text{min}}^{\text{(eff)}}$.
As can be seen in Fig.~\ref{fig:different-duration}, the cutoff
moves to the lower frequencies as $T$ increases, the power spectral
density is flat for the lowest observable natural frequencies, $\frac{1}{T}<f\lesssim\frac{\gamma_{\text{max}}}{\gamma_{\theta}T}$.

\begin{figure}[ht]
\begin{centering}
\includegraphics[width=0.45\textwidth]{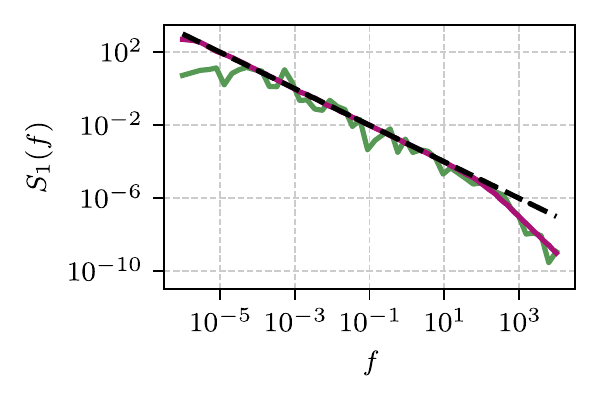}
\par\end{centering}
\caption{The effect of averaging over repeated experiments on the obtained
power spectral density: $R=1$ (green curve), $R=10^{3}$ (magenta
curve). Dashed black curve corresponds to Eq.~(\ref{eq:main-result}).
Simulation parameters, with exception to $R$, are the same as for
the green curve from Fig.~\ref{fig:different-duration}.\label{fig:different-duration-big-r}}
\end{figure}

If multiple independent experiments (let $R$ be the number of experiments)
with finite observation time $T$ are performed and the obtained spectral
densities are averaged, then the total number of observed pulses increases
by a factor of $R$ yielding
\begin{equation}
\gamma_{\text{min}}^{\text{(eff)}}=\left(\gamma_{\text{max}}-\gamma_{\text{min}}\right)\frac{\left\langle \theta\right\rangle +\left\langle \tau\right\rangle }{\left\langle \theta\right\rangle +\left\langle \tau\right\rangle +RT}+\gamma_{\text{min}}\approx\frac{\gamma_{\text{max}}\left\langle \theta\right\rangle }{RT}+\frac{1}{T}=\frac{R+\gamma_{\text{max}}\left\langle \theta\right\rangle }{RT}.
\end{equation}
For $R\gg\gamma_{\text{max}}\left\langle \theta\right\rangle $, no
low-frequency cutoff will be noticeable. As shown in Fig.~\ref{fig:different-duration-big-r},
low-frequency cutoff disappears as the experiments are repeated and
the obtained power spectral densities are averaged.

We have derived Eq.~(\ref{eq:main-result}) considering the current
generated by a single charge carrier. In many experiments the number
of charge carriers $N$ will be large, $N\gg1$. Consequently, from
the Wiener--Khinchin theorem \cite{Kogan1996CUP} it follows that
performing independent experiments is equivalent to observing independent
charge carriers. Therefore for $N\gg\gamma_{\text{max}}\left\langle \theta\right\rangle $
no low-frequency cutoff will be noticeable. Though in this case, the
power spectral densities of the signals generated by single charge
carriers add up instead of averaging out, yielding a minor generalization
of Eq.~(\ref{eq:main-result})
\begin{equation}
S_{N}\left(f\right)\approx\frac{Na^{2}\bar{\nu}}{\gamma_{\text{max}}f}.\label{eq:gen-main-result}
\end{equation}
In the above $\bar{\nu}$ is strictly the mean number of pulses per
unit time generated by a single charge carrier.

\begin{figure}[ht]
\begin{centering}
\includegraphics[width=0.9\textwidth]{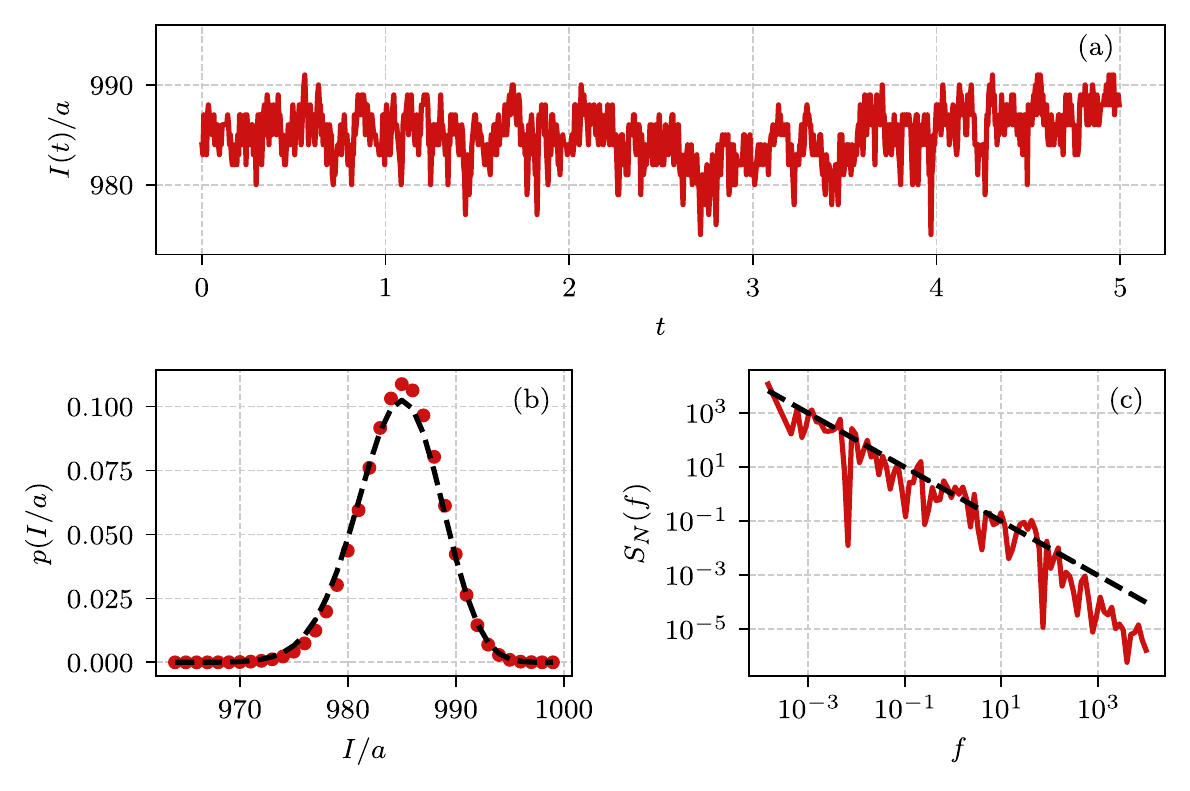}
\par\end{centering}
\caption{Results of a single simulation with large number of
charge carriers $N$ and finite duration $T$: excerpt of a signal
generated by $10^{3}$ independent charge carriers (a), the probability
mass function of the amplitude of the signal (b), and the power spectral
density of the signal (c). Red curves represent results of numerical
simulation, while dashed black curves provide theoretical fits: (b)
Binomial probability mass function with $p_{F}\approx0.984$ and $N=10^{3}$,
(c) the power spectral density approximation Eq.~(\ref{eq:gen-main-result}).
Simulation parameters: $R=1$, $N=10^{3}$, $T=2^{26}\cdot10^{-4}=6710.8864$,
$a=1$, $\gamma_{\theta}=1$, $\gamma_{\text{min}}=0$, $\gamma_{\text{max}}=10^{3}$.\label{fig:multicarrier}}
\end{figure}

As can be seen in Fig.~\ref{fig:multicarrier}~(a), the signal generated
by multiple independent charge carriers is no longer composed of non-overlapping
pulses, although it retains discrete nature as individual charges
drift freely or are trapped by the trapping centers. The amplitude
and the slope of the power spectral density are well predicted by
Eq.~(\ref{eq:gen-main-result}) (as seen in Fig.~\ref{fig:multicarrier}~(c)).
The distribution of the signal's amplitude would be expected to follow
the Binomial distribution with sample size $N$ and success probability
(probability that the charge carrier is free)
\begin{equation}
p_{F}=\frac{\left\langle \theta\right\rangle }{\left\langle \theta\right\rangle +\left\langle \tau\right\rangle }\approx1-\frac{\left\langle \tau\right\rangle }{\left\langle \theta\right\rangle }.
\end{equation}
The fit by the Binomial distribution shown in Fig.~\ref{fig:multicarrier}~(b)
is not perfect, because the nonergodic case is simulated and $\left\langle \tau\right\rangle $
is ill-defined, but predicts the overall shape of the probability
distribution rather well. For $\gamma_{\text{min}}\gg1/T$ the fit
would be much better. Notably, with larger $N$ and under noisy observation,
the Binomial distribution predicted by the model will quickly become
indistinguishable from the Gaussian distribution. While in some cases
$1/f$ noise is known to behave as a non-Gaussian process, most often
it is found to exhibit Gaussian fluctuations \cite{Kogan1996CUP,Melkonyan2010PhysB,Ruseckas2016JStatMech}.
The duration of the reported simulation was chosen arbitrarily, based
on the technical considerations. Specifically, we have opted to make
$2^{26}$ observations of the process with sampling period of $\Delta t=10^{-4}$.

Notably, \cite{Leibovich2017PRE} also discusses a spurious low-frequency
cutoff that could be observed in single particle experiments. Of the
$1/f$ noise models considered in \cite{Leibovich2017PRE} superimposed
random telegraph signals and blinking quantum dot models are the most
comparable to the model presented here. In \cite{Leibovich2017PRE}
each of the superimposed random telegraph signals was assumed to be
characterized by their own Poissonian switching rate $\gamma=\gamma_{\theta}=\gamma_{\tau}$
between the ``on'' and ``off'' states. It was shown that the conditional
power spectral density (requiring a certain minimum number of pulses,
$K_{\text{min}}$, to be observed) exhibits low-frequency cutoff at
$f_{c}\sim K_{\text{min}}/T$. In our simulations, we typically observe
a large number of pulses, $K\approx\gamma_{\theta}T$, and should
therefore observe the cutoff at $f_{c}\sim\gamma_{\theta}$, but instead,
we observe that the cutoff frequency scales as $1/\gamma_{\theta}$.
The nature of the cutoff is different in the model introduced here.
The other, blinking quantum dot, model does not predict low-frequency
cutoff, only the ageing effect, which for the pure $1/f$ noise will
not be noticeable \cite{Kononovicius2023PRE}.

\section{Derivation of Hooge's empirical relation and Hooge's parameter\label{sec:hooge-expression}}

It is straightforward to see that we can rewrite Eq.~(\ref{eq:gen-main-result})
in the form of Hooge's empirical relation, Eq.~(\ref{eq:hooge-relation}),
if we define Hooge's parameter as
\begin{equation}
\alpha_{H}=\frac{N^{2}a^{2}\bar{\nu}}{\gamma_{\text{max}}\bar{I}^{2}}.\label{eq:hooge-comined}
\end{equation}
Further we show that the straightforward expression above can be simplified,
and given a more compact form.

As the height of the pulses $a$ corresponds to the current generated
by a single charge carrier, we have
\begin{equation}
a=\frac{qv_{c}}{L},
\end{equation}
where $q$ stands for the charge held by the carrier, $v_{c}$ is
the free drift velocity between the trappings (which will be much
smaller than the thermal velocity of the charge carriers), and $L$
is the length of the material. Expression for $a$ can be rewritten
in terms of the average current flowing through the cross-section
of the material $\sigma_{M}$
\begin{equation}
\bar{I}=\sigma_{M}nqv_{d},
\end{equation}
where $n$ stands for the density of the charge carriers (i.e., $n=\frac{N}{L\sigma_{M}}$),
and $v_{d}$ is the average drift velocity of the charge carriers.
The average drift velocity is related to the free drift velocity via
the fraction of time the charge carrier spends drifting 
\begin{equation}
v_{d}=\frac{\left\langle \theta\right\rangle }{\left\langle \theta\right\rangle +\left\langle \tau\right\rangle }v_{c}=\bar{\nu}\left\langle \theta\right\rangle v_{c}.
\end{equation}
Consequently we have 
\begin{equation}
a=\frac{\bar{I}}{N\bar{\nu}\left\langle \theta\right\rangle }.\label{eq:hooge-pulse-height}
\end{equation}

Inserting Eq.~(\ref{eq:hooge-pulse-height}) into Eq.~(\ref{eq:hooge-comined})
yields the expression of the Hooge's parameter in terms of the characteristic
trapping rate and the maximum detrapping rate, assuming that the trapping
times are comparatively long $\left\langle \theta\right\rangle \gg\left\langle \tau\right\rangle $,
\begin{equation}
\alpha_{H}=\frac{1}{\bar{\nu}\left\langle \theta\right\rangle ^{2}\gamma_{\text{max}}}\approx\frac{\gamma_{\theta}}{\gamma_{\text{max}}}=\frac{\left\langle \tau_{\text{min}}\right\rangle }{\left\langle \theta\right\rangle }.\label{eq:hooge-final}
\end{equation}
In the above $\left\langle \tau_{\text{min}}\right\rangle =\frac{1}{\gamma_{\text{max}}}$
is the expected detrapping time generated when a charge carrier is
trapped by the shallowest trapping center. The purer materials (i.e.,
ones with lower trapping center density $n_{c}$) will have lower
$\alpha_{H}$ values, as the trapping rate is given by $\gamma_{\theta}=\left\langle \sigma_{c}v_{t}\right\rangle n_{c}$
(here $v_{t}$ is the thermal velocity of the charge carriers, and
$\sigma_{c}$ is the trapping cross-section). The proportionality
$\alpha_{H}\propto\gamma_{\theta}$ was previously reported in \cite{Lukyanchikova1990PhysB,Tacano2004FluctuationsNoiseMater,Tousek2024SciRep},
providing experimental support to Eq.~(\ref{eq:hooge-final}).

Consequently the approximations for the power spectral density generated
by the proposed model, Eqs.~(\ref{eq:main-result}) and (\ref{eq:gen-main-result}),
can be rewritten in the same form as Hooge's empirical relation. Inserting
Eq.~(\ref{eq:hooge-final}) into Eq.~(\ref{eq:hooge-relation})
yields
\begin{equation}
S_{N}\left(f\right)=\bar{I}^{2}\frac{\gamma_{\theta}}{\gamma_{\text{max}}Nf}.\label{eq:hooge-relation-final}
\end{equation}
This expression appears to imply that the process under consideration
is stationary, but this is not true as the average current $\bar{I}$
is proportional to the number of pulses per unit time $\bar{\nu}$,
which in the $\gamma_{\text{min}}\rightarrow0$ limit is a function
of the observation time $T$ \cite{Kononovicius2023PRE}. Although,
for the case of pure $1/f$ noise, the dependence on $T$ is logarithmically
slow, and barely noticeable. Nevertheless, even if the process would
be non-stationary, this should not have any impact on the estimate
of Hooge's parameter as only $\bar{I}$ is impacted by the non-stationarity.

\section{Conclusions\label{sec:conclusions}}

We have proposed a general model of $1/f$ noise in homogeneous semiconductors
which is based on the trapping-detrapping process of individual charge
carriers. In contrast to the many previous works, we have assumed
that the detrapping rate of each trapping center is random. We have
shown that, if detrapping process obeys Arrhenius law (which is well-established
empirically \cite{Peters2015JPCB,Kurpiers2018NatComm,Arakawa2020NatMater,Kumar2024PRE}),
and if the vacant trap depths follow Boltzman distribution (which
is also supported by experimental works \cite{Bisquert2008PRB,Wong2020ACSEnergyLett,Beckers2021JAppPhys,Zeiske2021NatureComm}),
the detrapping rate distribution will be uniform. When detrapping
rates are uniformly distributed, a power-law distribution of the detrapping
times Eq.~(\ref{eq:escape-time-pdf}) is obtained. It arises from
the superposition of exponential detrapping time distributions representing
contributions of the individual trapping centers with their own fixed
detrapping rates (see Fig\@.~\ref{fig:distribution}).

Consequently, regardless of the exact details of the trapping process,
as long as the trapping process is slow in comparison to the detrapping
process, pure $1/f$ noise in a form of Hooge's empirical relation
is obtained, Eq.~(\ref{eq:hooge-relation-final}). Corresponding
expression of the Hooge's parameter, $\alpha_{H}$, is then found
to be a ratio between the rate parameters of the trapping and the
detrapping processes, Eq.~(\ref{eq:hooge-final}). The proportionality
between the Hooge's parameter and the trapping rate was reported in
previous experimental works \cite{Lukyanchikova1990PhysB,Tacano2004FluctuationsNoiseMater,Tousek2024SciRep},
thus providing partial experimental verification for the Hooge's parameter
expression we have derived from general theoretical considerations.
Inverse proportionality between the Hooge's parameter and the maximum
detrapping rate suggests interesting implications for approaching
suppression of $1/f$ noise problem \cite{Kamada2023CommPhys,Nakatani2023JAppPhys,Balandin2024AppPhysLett}.
When the Arrhenius law applies, maximum detrapping rate could be increased
either by manipulating the pre-exponential factor, or by decreasing
shalowest trap depth (minimum activation energy) from which the Boltzmann
distribution applies to the trap depth distribution. The obtained
expression for the Hooge's parameter also suggests that $1/f$ noise
arises from the temporal charge carrier number fluctuations, not from
the spatial mobility fluctuations.

In Section~\ref{sec:finite-experiments}, we have discussed the implications
of finite experiments. We have shown that the power spectral density
may exhibit spurious low-frequency cutoff simply due to finite duration
of the experiment or simulation. The obtained width of the cutoff
is of the same order of magnitude as $\frac{\gamma_{\text{max}}}{\gamma_{\theta}}$.
This cutoff disappears when the power spectral density is averaged
over a large number of experiments, or when the experiment involves
a large number of independent charge carriers. In the latter case
the distribution of the signal's amplitude follows Binomial distribution,
which under imperfect observation will quickly become indistinguishable
from the Gaussian distribution.

All of the code used to perform the reported numerical simulations
is available at \cite{UpoissGithub}.

\begin{singlespace}
\section*{Author contributions}

\textbf{Aleksejus Kononovicius:} Software, Validation, Writing --
Original Draft, Writing -- Review \& Editing, Visualization. \textbf{Bronislovas
Kaulakys:} Conceptualization, Methodology, Writing -- Original Draft,
Writing -- Review \& Editing.


\begin{thebibliography}{10}
\expandafter\ifx\csname url\endcsname\relax
  \def\url#1{\texttt{#1}}\fi
\expandafter\ifx\csname urlprefix\endcsname\relax\def\urlprefix{URL }\fi
\expandafter\ifx\csname href\endcsname\relax
  \def\href#1#2{#2} \def\path#1{#1}\fi

\bibitem{Johnson1925PR}
J.~B. Johnson, The {S}chottky effect in low frequency circuits, Physical Review
  26 (1925) 71--85.
\newblock \href {https://doi.org/10.1103/PhysRev.26.71}
  {\path{doi:10.1103/PhysRev.26.71}}.

\bibitem{Schottky1926PR}
W.~Schottky, Small-shot effect and flicker effect, Physical Review 28~(1)
  (1926) 74--103.
\newblock \href {https://doi.org/10.1103/PhysRev.28.74}
  {\path{doi:10.1103/PhysRev.28.74}}.

\bibitem{Milotti2002}
E.~Milotti, 1/f noise: {A} pedagogical review (2002).
\newblock \href {http://arxiv.org/abs/physics/0204033}
  {\path{arXiv:physics/0204033}}, \href
  {https://doi.org/10.48550/arXiv.physics/0204033}
  {\path{doi:10.48550/arXiv.physics/0204033}}.

\bibitem{Kogan1996CUP}
S.~Kogan, Electronic noise and fluctuations in solids, Cambridge University
  Press, 1996.
\newblock \href {https://doi.org/10.1017/CBO9780511551666}
  {\path{doi:10.1017/CBO9780511551666}}.

\bibitem{Lowen2005Wiley}
S.~B. Lowen, M.~C. Teich, Fractal-Based Point Processes, Wiley, 2005.
\newblock \href {https://doi.org/10.1002/0471754722}
  {\path{doi:10.1002/0471754722}}.

\bibitem{VanKampen2007NorthHolland}
N.~G. van Kampen, Stochastic process in physics and chemistry, North Holland,
  Amsterdam, 2007.

\bibitem{Voss1976PRL}
R.~F. Voss, J.~Clarke, 1/f noise from systems in thermal equilibrium, Physical
  Review Letters 36~(1) (1976) 42--45.
\newblock \href {https://doi.org/10.1103/PhysRevLett.36.42}
  {\path{doi:10.1103/PhysRevLett.36.42}}.

\bibitem{Dutta1981RMP}
P.~Dutta, P.~M. Horn, Low-frequency fluctuations in solids: 1/f noise, Reviews
  of Modern Physics 53 (1981) 497--516.
\newblock \href {https://doi.org/10.1103/RevModPhys.53.497}
  {\path{doi:10.1103/RevModPhys.53.497}}.

\bibitem{Mitin2002}
V.~Mitin, L.~Reggiani, L.~Varani, Generation-recombination noise in
  semiconductors, in: Noise and Fluctuation Controls in Electronic Devices,
  Noise and Fluctuation Controls in Electronic Devices, American Scientific
  Publishers, 2002.

\bibitem{Balandin2013NN}
A.~A. Balandin, Low-frequency 1/f noise in graphene devices, Nature
  Nanotechnology 8 (2013) 549--555.
\newblock \href {https://doi.org/10.1038/nnano.2013.144}
  {\path{doi:10.1038/nnano.2013.144}}.

\bibitem{Fox2021Nature}
Z.~R. Fox, E.~Barkai, D.~Krapf, Aging power spectrum of membrane protein
  transport and other subordinated random walks, Nature Communications 12
  (2021).
\newblock \href {https://doi.org/10.1038/s41467-021-26465-8}
  {\path{doi:10.1038/s41467-021-26465-8}}.

\bibitem{Wirth2021IEEE}
G.~Wirth, M.~B. da~Silva, T.~H. Both, Unified compact modeling of charge
  trapping in 1/f noise, {RTN} and {BTI}, in: 2021 5th {IEEE} Electron Devices
  Technology and Manufacturing Conference ({EDTM}), {IEEE}, 2021, pp. 1--3.
\newblock \href {https://doi.org/10.1109/edtm50988.2021.9421005}
  {\path{doi:10.1109/edtm50988.2021.9421005}}.

\bibitem{Voss1975Nature}
R.~F. Voss, J.~Clarke, 1/f noise in music and speech, Nature 258 (1975)
  317--318.
\newblock \href {https://doi.org/10.1038/258317a0}
  {\path{doi:10.1038/258317a0}}.

\bibitem{Kobayashi1982BioMed}
M.~Kobayashi, T.~Musha, 1/f fluctuation of heartbeat period, IEEE Transactions
  on Biomedical Engineering 29 (1982) 456--457.
\newblock \href {https://doi.org/10.1109/TBME.1982.324972}
  {\path{doi:10.1109/TBME.1982.324972}}.

\bibitem{Gilden1995Science}
D.~L. Gilden, T.~Thornton, M.~W. Mallon, $1/ f$ noise in human cognition,
  Science 267~(5205) (1995) 1837--1839.
\newblock \href {https://doi.org/10.1126/science.7892611}
  {\path{doi:10.1126/science.7892611}}.

\bibitem{Levitin2012PNAS}
D.~J. Levitin, P.~Chordia, V.~Menon, Musical rhythm spectra from {B}ach to
  {J}oplin obey a 1/f power law, Proceedings of the National Academy of
  Sciences of the United States of America 109 (2012) 3716--3720.
\newblock \href {https://doi.org/10.1073/pnas.1113828109}
  {\path{doi:10.1073/pnas.1113828109}}.

\bibitem{Bernamont1937ProcPhysSoc}
J.~Bernamont, Fluctuations in the resistance of thin films, Proceedings of the
  Physical Society 49~(4S) (1937) 138--139.
\newblock \href {https://doi.org/10.1088/0959-5309/49/4S/316}
  {\path{doi:10.1088/0959-5309/49/4S/316}}.

\bibitem{Surdin1951}
M.~Surdin, Une th{\'{e}}orie des fluctuations {\'{e}}lectriques dans les
  semi-conducteurs, Journal de Physique et le Radium 12~(8) (1951) 777--783.
\newblock \href {https://doi.org/10.1051/jphysrad:01951001208077700}
  {\path{doi:10.1051/jphysrad:01951001208077700}}.

\bibitem{McWhorter1957}
A.~L. McWhorter, R.~H. Kingston, Semiconductor surface physics, in: Proceedings
  of the Conference on Physics of Semiconductor Surface Physics, Vol. 207,
  University of Pennsylvania, Philadelphia, 1957.

\bibitem{VanDerZiel1979AEEP}
A.~V.~D. Ziel, Flicker noise in electronic devices, in: Advances in Electronics
  and Electron Physics, Elsevier, 1979, pp. 225--297.
\newblock \href {https://doi.org/10.1016/s0065-2539(08)60768-4}
  {\path{doi:10.1016/s0065-2539(08)60768-4}}.

\bibitem{Kaulakys2005PhysRevE}
B.~Kaulakys, V.~Gontis, M.~Alaburda, Point process model of 1/f noise vs a sum
  of {L}orentzians, Physical Review E 71 (2005) 051105.
\newblock \href {https://doi.org/10.1103/PhysRevE.71.051105}
  {\path{doi:10.1103/PhysRevE.71.051105}}.

\bibitem{Wong2003MR}
H.~Wong, Low-frequency noise study in electron devices: {R}eview and update,
  Microelectronics Reliability 43~(4) (2003) 585--599.
\newblock \href {https://doi.org/10.1016/S0026-2714(02)00347-5}
  {\path{doi:10.1016/S0026-2714(02)00347-5}}.

\bibitem{Ruseckas2011EPL}
J.~Ruseckas, B.~Kaulakys, V.~Gontis, Herding model and 1/f noise, EPL 96 (2011)
  60007.
\newblock \href {https://doi.org/10.1209/0295-5075/96/60007}
  {\path{doi:10.1209/0295-5075/96/60007}}.

\bibitem{Kononovicius2012PhysA}
A.~Kononovicius, V.~Gontis, Agent based reasoning for the non-linear stochastic
  models of long-range memory, Physica A 391 (2012) 1309--1314.
\newblock \href {https://doi.org/10.1016/j.physa.2011.08.061}
  {\path{doi:10.1016/j.physa.2011.08.061}}.

\bibitem{Kaulakys2009JStatMech}
B.~Kaulakys, M.~Alaburda, Modeling scaled processes and $1/f^\beta$ noise using
  non-linear stochastic differential equations, Journal of Statistical
  Mechanics 2009~(02) (2009) P02051.
\newblock \href {https://doi.org/10.1088/1742-5468/2009/02/p02051}
  {\path{doi:10.1088/1742-5468/2009/02/p02051}}.

\bibitem{Kazakevicius2021Entropy}
R.~Kazakevicius, A.~Kononovicius, B.~Kaulakys, et~al., Understanding the nature
  of the long--range memory phenomenon in socioeconomic systems, Entropy 23
  (2021) 1125.
\newblock \href {https://doi.org/10.3390/e23091125}
  {\path{doi:10.3390/e23091125}}.

\bibitem{Ruseckas2014JStatMech}
J.~Ruseckas, B.~Kaulakys, Scaling properties of signals as origin of 1/f noise,
  Journal of Statistical Mechanics 2014~(6) (2014) P06005.
\newblock \href {https://doi.org/10.1088/1742-5468/2014/06/p06005}
  {\path{doi:10.1088/1742-5468/2014/06/p06005}}.

\bibitem{Kaulakys2015MPLB}
B.~Kaulakys, M.~Alaburda, J.~Ruseckas, 1/f noise from the nonlinear
  transformations of the variables, Modern Physics Letters B 29 (2015) 1550223.
\newblock \href {https://doi.org/10.1142/S0217984915502231}
  {\path{doi:10.1142/S0217984915502231}}.

\bibitem{Eliazar2021JPA}
I.~Eliazar, Selfsimilar diffusions, Journal of Physics A: Mathematical and
  Theoretical 54 (2021) 35LT01.
\newblock \href {https://doi.org/10.1088/1751-8121/ac1771}
  {\path{doi:10.1088/1751-8121/ac1771}}.

\bibitem{Kazakevicius2023PRE}
R.~Kazakevi{\v{c}}ius, A.~Kononovicius, Anomalous diffusion and long-range
  memory in the scaled voter model, Physical Review E 107 (2023) 024106.
\newblock \href {https://doi.org/10.1103/PhysRevE.107.024106}
  {\path{doi:10.1103/PhysRevE.107.024106}}.

\bibitem{Mandelbrot2013Springer}
B.~B. Mandelbrot, Multifractals and 1/f noise: Wild self-affinity in physics
  (1963--1976), Springer, 2013.

\bibitem{Beran2017Routledge}
J.~Beran, Statistics for long-memory processes, Routledge, 2017.
\newblock \href {https://doi.org/10.1201/9780203738481}
  {\path{doi:10.1201/9780203738481}}.

\bibitem{Giraitis2009}
L.~Giraitis, R.~Leipus, D.~Surgailis, {ARCH}($\infty$) models and long memory,
  in: T.~G. Anderson, R.~A. Davis, J.~Kreis, T.~Mikosh (Eds.), Handbook of
  Financial Time Series, Springer Verlag, Berlin, 2009, pp. 71--84.
\newblock \href {https://doi.org/10.1007/978-3-540-71297-8\_3}
  {\path{doi:10.1007/978-3-540-71297-8\_3}}.

\bibitem{Bollerslev2023JEconom}
T.~Bollerslev, The story of {GARCH}: {A} personal odyssey, Journal of
  Econometrics 234 (2023) 96--100.
\newblock \href {https://doi.org/10.1016/j.jeconom.2023.01.015}
  {\path{doi:10.1016/j.jeconom.2023.01.015}}.

\bibitem{Hooge1969PLA}
F.~N. Hooge, 1/f noise is no surface effect, Physics Letters A 29~(3) (1969)
  139--140.
\newblock \href {https://doi.org/10.1016/0375-9601(69)90076-0}
  {\path{doi:10.1016/0375-9601(69)90076-0}}.

\bibitem{Hooge1972Phys}
F.~N. Hooge, Discussion of recent experiments on 1/f noise, Physica 60~(1)
  (1972) 130--144.
\newblock \href {https://doi.org/10.1016/0031-8914(72)90226-1}
  {\path{doi:10.1016/0031-8914(72)90226-1}}.

\bibitem{Hooge1994IEEE}
F.~N. Hooge, 1/f noise sources, {IEEE} Transactions on Electron Devices 41~(11)
  (1994) 1926--1935.
\newblock \href {https://doi.org/10.1109/16.333808}
  {\path{doi:10.1109/16.333808}}.

\bibitem{Kaulakys1999PLA}
B.~Kaulakys, Autoregressive model of 1/f noise, Physics Letters A 257 (1999)
  37--42.
\newblock \href {https://doi.org/10.1016/S0375-9601(99)00284-4}
  {\path{doi:10.1016/S0375-9601(99)00284-4}}.

\bibitem{Dmitriev2009JAppPhys}
A.~P. Dmitriev, M.~E. Levinshtein, S.~L. Rumyantsev, On the hooge relation in
  semiconductors and metals, Journal of Applied Physics 106~(2) (2009) 024514.
\newblock \href {https://doi.org/10.1063/1.3186620}
  {\path{doi:10.1063/1.3186620}}.

\bibitem{Vandamme2013ICNF}
L.~K.~J. Vandamme, How useful is {H}ooge's empirical relation, in: 22nd
  International Conference on Noise and Fluctuations ({ICNF}), {IEEE}, 2013,
  pp. 1--6.
\newblock \href {https://doi.org/10.1109/ICNF.2013.6578875}
  {\path{doi:10.1109/ICNF.2013.6578875}}.

\bibitem{Palenskis2015SRep}
V.~Palenskis, K.~Maknys, Nature of low-frequency noise in homogeneous
  semiconductors, Scientific Reports 5~(1) (2015).
\newblock \href {https://doi.org/10.1038/srep18305}
  {\path{doi:10.1038/srep18305}}.

\bibitem{Gruneis2019PLA}
F.~Gruneis, An alternative form of {H}ooge's relation for 1/f noise in
  semiconductor materials, Physics Letters A 383~(13) (2019) 1401--1409.
\newblock \href {https://doi.org/10.1016/j.physleta.2019.02.009}
  {\path{doi:10.1016/j.physleta.2019.02.009}}.

\bibitem{Gruneis2022PhysA}
F.~Gruneis, 1/f noise under drift and thermal agitation in semiconductor
  materials, Physica A 593 (2022) 126917.
\newblock \href {https://doi.org/10.1016/j.physa.2022.126917}
  {\path{doi:10.1016/j.physa.2022.126917}}.

\bibitem{Bisquert2008PRB}
J.~Bisquert, Beyond the quasistatic approximation: {I}mpedance and capacitance
  of an exponential distribution of traps, Physical Review B 77~(23) (2008)
  235203.
\newblock \href {https://doi.org/10.1103/PhysRevB.77.235203}
  {\path{doi:10.1103/PhysRevB.77.235203}}.

\bibitem{Wong2020ACSEnergyLett}
J.~Wong, S.~T. Omelchenko, H.~A. Atwater, Impact of semiconductor band tails
  and band filling on photovoltaic efficiency limits, ACS Energy Letters 6~(1)
  (2020) 52--57.
\newblock \href {https://doi.org/10.1021/acsenergylett.0c02362}
  {\path{doi:10.1021/acsenergylett.0c02362}}.

\bibitem{Beckers2021JAppPhys}
A.~Beckers, D.~Beckers, F.~Jazaeri, et~al., Generalized {B}oltzmann relations
  in semiconductors including band tails, Journal of Applied Physics 129~(4)
  (2021).
\newblock \href {https://doi.org/10.1063/5.0037432}
  {\path{doi:10.1063/5.0037432}}.

\bibitem{Zeiske2021NatureComm}
S.~Zeiske, O.~J. Sandberg, N.~Zarrabi, et~al., Direct observation of
  trap-assisted recombination in organic photovoltaic devices, Nature
  Communications 12~(1) (2021).
\newblock \href {https://doi.org/10.1038/s41467-021-23870-x}
  {\path{doi:10.1038/s41467-021-23870-x}}.

\bibitem{Peters2015JPCB}
B.~Peters, Common features of extraordinary rate theories, The Journal of
  Physical Chemistry B 119~(21) (2015) 6349--6356.
\newblock \href {https://doi.org/10.1021/acs.jpcb.5b02547}
  {\path{doi:10.1021/acs.jpcb.5b02547}}.

\bibitem{Kurpiers2018NatComm}
J.~Kurpiers, T.~Ferron, S.~Roland, et~al., Probing the pathways of free charge
  generation in organic bulk heterojunction solar cells, Nature Communications
  9~(1) (2018) 2038.
\newblock \href {https://doi.org/10.1038/s41467-018-04386-3}
  {\path{doi:10.1038/s41467-018-04386-3}}.

\bibitem{Arakawa2020NatMater}
K.~Arakawa, M.-C. Marinica, S.~Fitzgerald, et~al., Quantum de-trapping and
  transport of heavy defects in tungsten, Nature Materials 19~(5) (2020)
  508--511.
\newblock \href {https://doi.org/10.1038/s41563-019-0584-0}
  {\path{doi:10.1038/s41563-019-0584-0}}.

\bibitem{Kumar2024PRE}
V.~Kumar, A.~Pal, O.~Shpielberg, Arrhenius law for interacting diffusive
  systems, Physical Review E 109~(3) (2024) l032101.
\newblock \href {https://doi.org/10.1103/PhysRevE.109.L032101}
  {\path{doi:10.1103/PhysRevE.109.L032101}}.

\bibitem{Kononovicius2023PRE}
A.~Kononovicius, B.~Kaulakys, 1/f noise from the sequence of nonoverlapping
  rectangular pulses, Physical Review E 107 (2023) 034117.
\newblock \href {https://doi.org/10.1103/PhysRevE.107.034117}
  {\path{doi:10.1103/PhysRevE.107.034117}}.

\bibitem{Lukyanchikova1990PhysB}
N.~Lukyanchikova, M.~Petrichuk, N.~Garbar, et~al., 1/f noise and
  generation-recombination processes at discrete levels in semiconductors,
  Physica B: Condensed Matter 167~(3) (1990) 201--207.
\newblock \href {https://doi.org/10.1016/0921-4526(90)90352-u}
  {\path{doi:10.1016/0921-4526(90)90352-u}}.

\bibitem{Tacano2004FluctuationsNoiseMater}
M.~Tacano, J.~Pavelka, N.~Tanuma, S.~Yokokura, S.~Hashiguchi, Dependence of
  {H}ooge constant on mean free paths of materials, in: D.~Popovic, M.~B.
  Weissman, Z.~A. Racz (Eds.), Fluctuations and Noise in Materials, SPIE, 2004,
  pp. 310--319.
\newblock \href {https://doi.org/10.1117/12.547202}
  {\path{doi:10.1117/12.547202}}.

\bibitem{Tousek2024SciRep}
J.~Tousek, J.~Touskova, I.~Krivka, Product of mobility and lifetime of charge
  carriers in cdte determined from low-frequency current fluctuations,
  Scientific Reports 14~(1) (2024).
\newblock \href {https://doi.org/10.1038/s41598-024-51541-6}
  {\path{doi:10.1038/s41598-024-51541-6}}.

\bibitem{Margolin2006JStatPhys}
G.~Margolin, E.~Barkai, Nonergodicity of a time series obeying {L}evy
  statistics, Journal of Statistical Physics 122~(1) (2006) 137--167.
\newblock \href {https://doi.org/10.1007/s10955-005-8076-9}
  {\path{doi:10.1007/s10955-005-8076-9}}.

\bibitem{Lukovic2008JChemPhys}
M.~Lukovic, P.~Grigolini, Power spectra for both interrupted and perennial
  aging processes, The Journal of Chemical Physics 129~(18) (2008) 184102.
\newblock \href {https://doi.org/10.1063/1.3006051}
  {\path{doi:10.1063/1.3006051}}.

\bibitem{Niemann2013PRL}
M.~Niemann, H.~Kantz, E.~Barkai, Fluctuations of 1/f noise and the
  low-frequency cutoff paradox, Physical Review Letters 110~(14) (2013) 140603.
\newblock \href {https://doi.org/10.1103/PhysRevLett.110.140603}
  {\path{doi:10.1103/PhysRevLett.110.140603}}.

\bibitem{Leibovich2016PRE}
N.~Leibovich, A.~Dechant, E.~Lutz, et~al., Aging {W}iener-{K}hinchin theorem
  and critical exponents of $1/f^\beta$ noise, Physical Review E 94 (2016)
  052130.
\newblock \href {https://doi.org/10.1103/PhysRevE.94.052130}
  {\path{doi:10.1103/PhysRevE.94.052130}}.

\bibitem{Frantsuzov2008NatPhys}
P.~Frantsuzov, M.~Kuno, B.~Jank{\'{o}}, et~al., Universal emission
  intermittency in quantum dots, nanorods and nanowires, Nature Physics 4~(7)
  (2008) 519--522.
\newblock \href {https://doi.org/10.1038/nphys1001}
  {\path{doi:10.1038/nphys1001}}.

\bibitem{Cordones2013CSR}
A.~A. Cordones, S.~R. Leone, Mechanisms for charge trapping in single
  semiconductor nanocrystals probed by fluorescence blinking, Chemical Society
  Reviews 42~(8) (2013) 3209.
\newblock \href {https://doi.org/10.1039/c2cs35452g}
  {\path{doi:10.1039/c2cs35452g}}.

\bibitem{Nenashev2018PRB}
A.~V. Nenashev, V.~V. Valkovskii, J.~O. Oelerich, et~al., Release of carriers
  from traps enhanced by hopping, Physical Review B 98~(15) (2018) 155207.
\newblock \href {https://doi.org/10.1103/PhysRevB.98.155207}
  {\path{doi:10.1103/PhysRevB.98.155207}}.

\bibitem{Haneef2020JMCC}
H.~F. Haneef, A.~M. Zeidell, O.~D. Jurchescu, Charge carrier traps in organic
  semiconductors: {A} review on the underlying physics and impact on electronic
  devices, Journal of Materials Chemistry C 8~(3) (2020) 759--787.
\newblock \href {https://doi.org/10.1039/c9tc05695e}
  {\path{doi:10.1039/c9tc05695e}}.

\bibitem{Gentle2009Springer}
J.~E. Gentle, Mathematical and statistical preliminaries, in: Statistics and
  Computing, Springer New York, 2009, pp. 5--79.
\newblock \href {https://doi.org/10.1007/978-0-387-98144-4\_1}
  {\path{doi:10.1007/978-0-387-98144-4\_1}}.

\bibitem{Melkonyan2010PhysB}
S.~V. Melkonyan, Non-{G}aussian conductivity fluctuations in semiconductors,
  Physica B: Condensed Matter 405~(1) (2010) 379--385.
\newblock \href {https://doi.org/10.1016/j.physb.2009.08.096}
  {\path{doi:10.1016/j.physb.2009.08.096}}.

\bibitem{Ruseckas2016JStatMech}
J.~Ruseckas, R.~Kazakevicius, B.~Kaulakys, Coupled nonlinear stochastic
  differential equations generating arbitrary distributed observable with 1/f
  noise, Journal of Statistical Mechanics 2016~(4) (2016) 043209.
\newblock \href {https://doi.org/10.1088/1742-5468/2016/04/043209}
  {\path{doi:10.1088/1742-5468/2016/04/043209}}.

\bibitem{Leibovich2017PRE}
N.~Leibovich, E.~Barkai, Conditional $1/f^\alpha$ noise: {F}rom single
  molecules to macroscopic measurement, Physical Review E 96 (2017) 032132.
\newblock \href {https://doi.org/10.1103/PhysRevE.96.032132}
  {\path{doi:10.1103/PhysRevE.96.032132}}.

\bibitem{Kamada2023CommPhys}
M.~Kamada, W.~Zeng, A.~Laitinen, et~al., Suppression of 1/f noise in graphene
  due to anisotropic mobility fluctuations induced by impurity motion,
  Communications Physics 6~(1) (2023).
\newblock \href {https://doi.org/10.1038/s42005-023-01321-x}
  {\path{doi:10.1038/s42005-023-01321-x}}.

\bibitem{Nakatani2023JAppPhys}
T.~Nakatani, H.~Suto, P.~D. Kulkarni, et~al., Improvement of magnetic field
  detectivity in electrical 1/f noise-dominated tunnel magnetoresistive sensors
  by {AC} magnetic field modulation technique, Journal of Applied Physics
  134~(21) (2023).
\newblock \href {https://doi.org/10.1063/5.0180812}
  {\path{doi:10.1063/5.0180812}}.

\bibitem{Balandin2024AppPhysLett}
A.~A. Balandin, E.~Paladino, P.~J. Hakonen, Electronic noise---{F}rom advanced
  materials to quantum technologies, Applied Physics Letters 124~(5) (2024).
\newblock \href {https://doi.org/10.1063/5.0197142}
  {\path{doi:10.1063/5.0197142}}.

\bibitem{UpoissGithub}
\urlprefix\url{https://github.com/akononovicius/flicker-trap-detrap-individual-charge}

\end{thebibliography}

\end{singlespace}

\end{document}